\documentclass{amsart}
\usepackage{amssymb}

\newcommand{\R}{{\mathbb R}}

\newtheorem{definition}{Definition}[section]

\newtheorem{theorem}[definition]{Theorem}

\theoremstyle{remark}

\numberwithin{equation}{section}

\makeatletter
\begin{document}
\title{Parabolic equations with rough data}
\author{Herbert Koch}
\address[H.~Koch]{Mathematisches Institut\\Universit{\"a}t Bonn\\Endenicher Allee 60\\53115 Bonn\\Germany}
\email{koch@math.uni-bonn.de}
\author{Tobias Lamm}
\address[T.~Lamm]{Institute for Analysis, Karlsruhe Institute of Technology (KIT), Kaiserstr. 89-93, D-76133 Karlsruhe, Germany}
\email{tobias.lamm@kit.edu}

\begin{abstract}
We survey recent work on local well-posedness results for parabolic equations and systems with rough initial data.
\end{abstract}
\maketitle

\section{Introduction}\label{sec1}

In this paper we survey recent work on the initial value problem
for parabolic equations in a fairly broad sense. This new approach is based
on basic notions in harmonic analysis like maximal function, square function,
and Carleson measures. The design of the function spaces we use is modeled
on maximal functions and square functions, where the version we use incorporates
the regularity theory for the corresponding linear parabolic equations.

We consider it to be an appealing feature that a first local existence
statement can be formulated without using function spaces, while being
essentially optimal in terms of the regularity of the initial data
needed, see Theorem \ref{wellposed} below.

Our proofs make only use of fairly general properties of linear
equations with constant coefficients: (Gaussian) decay of the kernel,
and a version of the Calderon-Zygmund estimates.  Moreover the
arguments are almost local in space for local in time solutions. In
the flat small data situation this idea has first been used in Koch
and Tataru \cite{MR1808843} and, closer to the chore of this survey,
by the authors in \cite{MR2916362}.

One of the main observation is that the methods we use are flexible enough to handle are initial boundary value problems
in half spaces, parabolic systems, subelliptic parabolic equations,
and higher order parabolic equations.  Subelliptic parabolic equations
occur in the context of the porous medium equation (see \cite{zbMATH01198486} and the thesis of C. Kienzler \cite{kienzler}) and in the context of thin films (see the thesis of D. John \cite{john}).

It seems natural to study parabolic equations on uniform manifolds - i.e.
manifolds with a metric and an atlas corresponding to balls of size one
for which all the coordinate changes are uniformly in $C^1$ with uniform
modulus of continuity. 
This concept of uniform manifolds has been introduced by Denzler, Koch and McCann \cite{arxiv12046434} and it was recently used by  Shao and Simonett \cite{arxiv13092041} and Shao \cite{arxiv13092042}.

It is a consequence of our results- and basically this result can also be found in the papers of Whitney \cite{MR1503315} and Kotschwar \cite{zbMATH06137850} -
 that those manifolds carry a uniform
analytic metric: there is an atlas corresponding to balls of diameter
$1$ and a metric $g$ so that all coordinate changes $\phi_{ij}$ satisfy bounds
\[  | \partial_x^\alpha \phi_{ij}| \le  c R^{-|\alpha|}|\alpha|! \]
\[ |\partial_x^\alpha g^{ij} |  \le c R^{-|\alpha|} |\alpha|! \]
where $c$ and $R$ are independent of $\alpha$.

Initial boundary value problems fit into the 
framework of uniform structures:
Consider a bounded domain with smooth boundary.
Locally we can flatten the boundary, and we obtain a 'uniform' structure 
in the spirit as discussed above.

In the following we discuss several examples which we consider
instructive and interesting. Details will appear in \cite{kl13}.

Consider the equation
\begin{align} 
u_t - \sum_{i,j=1}^d \partial_i a^{ij}(t,x,u) \partial_j u = f(t,x,u,\nabla u) \label{zero}
\end{align}
in $\R^d$ where $a^{ij}$ and $f$ are continuous functions satisfying
\[  \lambda^{-1} |\xi|^2 \le \sum_{i,j=1}^d a^{ij}(t,x,u) \xi_i\xi_j \]
and
\[  |a^{ij}| \le \lambda \]
for some $\lambda>1$, uniformly for all $t,x,u$ and $\xi$. The
coefficients are not assumed to be symmetric.

The basic regularity assumption with respect to $x$ and $t$ is the requirement of
{\em locally small oscillation:}
There exists $\delta$ depending on $\lambda$,  and $T>0$ with
\[  |a^{ij}(t,x,u) -a^{ij}(s,y,u)| \le \delta \ \ \forall \ \ 0\le s,t \le T,
|x-y| \le \sqrt{T}  \]

We assume {\em Lipschitz continuity} with respect to $u$:
There exists $L$ with
\[ |a^{ij}(t,x,u)-a^{ij}(t,x,v)| \le L |u-v|. \]
The nonlinearity $f$ is assumed to be quadratic in the last component.
There is a small parameter $\varepsilon$ and we assume
\[ |f(t,x,u,0)| \le  \varepsilon/T \]
and
 \[ |f(t,x,u,p)-f(t,x,v,q)| \le  c\Big( |u-v| ( \varepsilon/T +L |p|^2)
+ (\varepsilon/\sqrt{T} + L(|p|+|q|))|p-q| \Big). \]

{\em Higher regularity:} Let $k\ge 1$ be a regularity index.
The  derivatives of $a^{ij}$ with respect
to $x$ and $u$ of order $k$ are uniformly bounded:
\[  T^{|\alpha|/2}  |\partial_x^{\alpha} \partial_u^{j}  a^{ij} | \le L  \]
and
\[  T^{1+|\alpha|/2-|\beta|/2} |\partial_x^{\alpha} \partial_u^{j} \partial_p^{\beta} f| \le L
 (1 + |T^{\frac12} p|^{(2-|\beta|)_+}).  \]
for $|\alpha|+j + |\beta| \le k$.

\begin{theorem}\label{wellposed}  There exists $\delta>0$, and for all $L>0$ there is $\varepsilon_0>0$ so that, if for $T>0$ 
\[
|u_0(x)-u_0(y)| \le \varepsilon  < \varepsilon_0 \text{ for  } |x-y| \le \sqrt{T}
\]
and the assumptions above are satisfied 
then there is a unique continuous solution $u$ up to time $T$ which satisfies
\[  |(t^{\frac12} \partial_x)^\alpha u(t,x) | \le c_{\alpha}  \varepsilon \]
for $|\alpha| \le k$. The solution is analytic with respect
to $x$ if $a^{ij}$  and $f$ are analytic.
If $a^{ij}$ and $f$ are analytic with respect to all variables then there exist $c$ and $R$ so that
\[  |(t^{\frac12}  \partial_x)^\alpha (t\partial_t)^j u(t,x) | \le c (|\alpha|+j)! R^{|\alpha|+j}  \varepsilon. \]
\end{theorem}

Examples of equations and systems of the above type are the harmonic
map heat flow, the viscous Hamilton Jacobi equation, the Ricci-DeTurck
flow and the fast diffusion equations for the relative size with
respect to the Barenblatt solution. In all of these cases continuous
initial data are natural and essentially optimal, which can be seen by
the examples below.

\section{The fixed point formulation}
\label{section2}
We construct the solution of the parabolic equation as a fixed point using Duhamel's formula.
For this we consider the abstract equation
\[ u_t = Au + f[u] \]
where $A$ is the generator of a semigroup $S(t)$.
If there are function spaces $X_0$, $X$ and $Y$ so that
\begin{equation}    \Vert S(t) u_0 \Vert_{X} \le c \Vert u_0 \Vert_{X_0} \end{equation}
\begin{equation}    \left\Vert \int_0^t  S(t-s) f(s) ds \right\Vert_X \le c \Vert f \Vert_Y \end{equation}
\begin{equation}  \Vert f[u]-f[v] \Vert_Y \le c (\Vert u \Vert_X + \Vert v \Vert_X+ \delta)
\Vert u-v \Vert_{X}, \end{equation}
then it is standard to deduce
\begin{itemize}
\item Existence and uniqueness by the contraction mapping principle.
\item Alternatively existence of the fixed point follows from the
  implicit function theorem, provided the maps are differentiable.  The
  contraction property implies invertibility of the
  linearization. This has an important consequence: The solution
  depends smoothly on parameters - if the nonlinear functions are
  smooth - resp. analytic if the functions are analytic.
\end{itemize}

Possible and popular  choices are
\begin{itemize}
\item H\"older spaces $C^\alpha(\Omega)$ and $C^{\alpha/2,\alpha}([0,T) \times \Omega)$ (see \cite{arxiv13092041}  for a recent contribution, discussion and references) 
\item The Sobolev space $X= W^{1,2,p}([0,T)\times \Omega)$, $X_0 =
  W^{2-\frac2p,p}(\Omega)$ of functions with one time and two spacial
  derivatives in $L^2$, $Y = L^p$, $p>n+2$
\end{itemize}

To motivate our choice we take a look at fundamental objects in
harmonic analysis. Consider the heat equation
\[ u_t = \Delta u, \quad u(0,x) = v(x) . \]
A nontangential maximal function is given by
\[ Mv(x) = \sup_{|h|^2\le t} |u(x+h,t)|   \]
which has the variant for $k \ge 0$ and $p\in [1,\infty]$
\[
 Mv(x) = \sup_{R}  R^k \left( R^{-d-2}    \int_{R^2/2}^{R^2}  \int_{B_{R}(x)}  |D^k_x  u|^p dx \, dt  \right)^{\frac1p}.
\]
The basic property is
\[   \Vert  Mv \Vert_{L^p} \le c \Vert v \Vert_{L^p} \]
for $1 < p \le \infty$ and
\[ \Vert v \Vert_{L^p} \le c \Vert Mv \Vert_{L^p} \]
if $1 <p < \infty$.

For $p= \infty$ there is a substitute via the square function
\[ \Vert v \Vert_{BMO} \sim
\sup_{x,R} \left(
R^{-d} \int_0^{R^2}   \int_{B_{R}(x)}   |\nabla u|^2 dy dt \right)^{\frac12}. \]
The right hand side is a Carleson measure type expression.

These tools have been used in the study of function spaces, but also
for the solution of the Kato square root problem by Auscher,
Hofmann, Lacey, McIntosh and Tchamitchian \cite{zbMATH01851146} and the study of harmonic functions
in Lipschitz domains by  Jerison, Kenig \cite{zbMATH00764042} and others.

The Carleson measure formulation of the $BMO$-norm (or more precisely the $BMO^{-1}$-norm) turned out to be a crucial ingredient in the study of the
 Navier-Stokes equations with initial data in $BMO^{-1}$ by Koch \& Tataru 
\cite{MR1808843} . More recently the authors applied these concepts to geometric problems including the harmonic map heat flow,
the Ricci-DeTurck flow, and the mean curvature and Willmore flow for Lipschitz graphs, see \cite{MR2916362}.

In order to study equations of the form \eqref{zero} we pick $p> n+2$ and $ q=p/2$. Moreover, we let $T>0$ and define the norms
\[  \Vert u_0 \Vert_{X_0}  = \Vert u_0 \Vert_{sup} \]
and
\[
\begin{split}
 \Vert u\Vert_{X} = & \sup_{x, t \le T} |u(t,x) |
\\ & +  \sup_{x, R^2< T} R  \left( R^{-d-2} \int_{R^2/2}^{R^2}\int_{B_R(x)}  |\nabla u |^p dy dt \right)^{\frac1p}
\\ &+ \sup_{x,R^2 < T }  \left( R^{-d} \int_{0}^{R^2}\int_{B_R(x)}  |\nabla u |^2 dy dt \right)^{\frac12}.
 \end{split}
\]
Here the second line is similar to the $L^\infty$ norm of a maximal function, 
and the last line corresponds to a Carleson measure. 

Additionally, we consider nonlinearities
\[  f[u] =  f_0(u,\nabla u) + \partial_i F^i(u,\nabla u)  \]
and decompose
\[ \Vert f \Vert_{Y} = \Vert f_0 \Vert_{Y^0} + \Vert F \Vert_{Y^1}, \]
where
\[ \begin{split} \Vert f_0 \Vert_{Y^0} =  &   \sup_{x, R^2< T} R \left(
R^{-\frac{d+2}2}\int_{R^2/2}^{R^2}\int_{B_R(x)}  |f_0|^q dy dt \right)^{\frac1q}
\\  &+ \sup_{x,R^2 < T }   R^{-d} \int_{0}^{R^2}\int_{B_R(x)}  |f_0| dy dt
\end{split}
\]
and
\[ \begin{split} 
\Vert F \Vert_{Y^1} = &
 \sup_{x, R^2< T}  R \left( R^{-d-2} \int_{R^2/2}^{R^2}\int_{B_R(x)}
|F|^p dy dt \right)^{\frac1p}\\
&+ \sup_{x,R^2 < T }  \left( R^{-d} \int_{0}^{R^2}\int_{B_R(x)}  |F |^2 dy dt \right)^{\frac12}.\end{split} 
\]

Now we construct a function $w:\R^d\to \R$ so that $||w-u_0||_{sup}$ and $||\nabla w||_{sup}$ are small in terms of $\varepsilon$, and we look for $u$ solving
\[
u_t - \sum_{i,j=1}^d \partial_i a^{ij}(t,x,w) \partial_j u
= f(t,x,u,\nabla u) + \sum_{i,j=1}^d \partial_i (a^{ij}(t,x,u) -a^{ij}(t,x,w)) \partial_j u)
\]
The estimate
\[ \Vert u \Vert_X \le c \Vert u_0 \Vert_{X_0} \]
follows from standard kernel estimates.
The estimates
\[
\Vert f(u)-f(v) \Vert_{Y^0} \le c (\Vert u \Vert_{X}+ \Vert v \Vert_X
+ \delta) \Vert u-v \Vert_{X}
\]
and the bound for $ \Vert
(a^{ij}(t,x,u) -a^{ij}(t,x,w)) \partial_j u \Vert_{Y_1} $ are true by
construction.

By scaling and the kernel estimates, if
\[  u_t - \sum_{i,j=1}^d \partial_i a^{ij} \partial_j u
= f + \sum_{i=1}^d \partial_i F^i
\]
with $u(0) = 0 $, then
\[ |u(0,1)| \le c \Big[\Vert f \Vert_{Y^0} + \Vert F \Vert_{Y^1}\Big]. \]
Energy estimates (plus kernel estimates) give
\[\left(  \int_0^1 \int_{B_1(0)}  |\nabla u|^2 dx dt  \right)^{\frac12}
\le c \Big[\Vert f \Vert_{Y^0} + \Vert F \Vert_{Y^1}\Big]. \]
Finally, kernel estimates and Calderon-Zygmund theory imply
\[ \left(  \int_{\frac12} ^1 \int_{B_1(0)}  |\nabla u|^p dx dt  \right)^{\frac1p}
\le c \Big[\Vert f \Vert_{Y^0} + \Vert F \Vert_{Y^1}\Big] . \]
Hence, we are in the above mentioned abstract setting in which the existence and uniqueness of a solution of \eqref{zero} in $X$ follows from a fixed point argument.

\section{Regularity and uniqueness}

We will prove regularity and existence of derivatives via the implicit
function theorem. For simplicity we do that for analyticity, where
this goes back to Angenent (see \cite{MR1078266}, \cite{MR1059647}).
Consider
\[ u_t - \Delta u = \Gamma(u) |\nabla u |^2 \]
for some analytic and bounded function $\Gamma$.
Define
\[ u^{s,a}(t,x)  = u(st, x+ta). \]
It satisfies
\[
u_t - s \Delta u + a \cdot \nabla u
= \Gamma(u) |\nabla u|^2
\]
which is analytic in $a$ (for $a$ close to zero) and $s$ (for $s$ close to $1$). We construct the solution by the implicit
function theorem. Thus $(s,a) \to u^{s,a} \in X $ is analytic.
The evaluation of a derivative is linear, hence for all $t$ and $x$ the map
\[ (s,a) \to \nabla u^{s,a}(t,x) \]
is analytic. But
\[ t\partial_t u = \partial_s u^{s,a}\Big|_{s=1,a=0} \]
and
\[ t \partial_j u = \partial_{a_j}  u^{s,a}\Big|_{s=1,a=0}, \]
with corresponding formulas for higher derivatives.

This argument can be localized as follows: The map
\[ x \to x+ ta \]
is the flow map of the constant vector field $a\in \R^d$. This is clearly
analytic with respect to $a$. We fix the analytic vector fields (for given $a$)
\[ X = (1-|x|^2)a \]
They generate a flow which is analytic with respect to $x$ and $a$.
The vector field vanishes at $|x|=1$. Hence also
\[ X_+ = (1-|x|^2)_+ a \]
generates a $C^{\infty}$ flow which is analytic with respect to the
parameter $a\in \R^d$. This argument shows that analyticity with
respect to $x$ is a local problem - in contrast to analyticity in
time: The fundamental solution is smooth but not analytic at $t=0$ and
$ x \in \R^d$ , $ x \ne 0$.

We need slightly more for Theorem \ref{wellposed}.  The properties on the nonlinearity are too
weak for a direct implementation of this argument. In order to overcome this difficulty we first obtain
bounds for the derivatives, and then we implement this argument on the
second half of the time interval.

Let us further comment on the uniqueness result claimed in Theorem \ref{wellposed}.
The fixed point map gives a unique fixed point in $X$ but Theorem \ref{wellposed}
claims uniqueness for weak solutions satisfying
\[t^{\frac12}  \Vert \nabla u \Vert_{L^{\infty}} \le c_1\varepsilon \]
which does not imply  the Carleson measure bound.  
Let $u$ be a solution as in the theorem.
For $t>0$ we can solve the initial value problem for the initial data
$u(t)$. It is unique, and hence the shifted solution is uniformly bounded
in $X$. The limit $t \to 0$ shows that the solution is in $X$ and hence unique.

\bigskip
The above framework allows to deal with rougher initial data.  It is obvious
that we may allow small perturbations of the initial data in
$L^\infty$. We may also allow small $BMO$ perturbations, if we require
that all the structure assumptions hold uniformly in $u$.  Here we
take a caloric extension $w$ of the initial data, and make the ansatz $u= w+v$. Then we apply a fixed point argument in order to find a function $v$ so that $u$ is a solution of our problem.

\section{Modifications and generalizations}

\subsection{Uniform manifolds and initial boundary value problems}

Parabolic equations have an infinite speed of propagation but heat
kernels have Gaussian decay.  Therefore we only need local in time estimates if
we want to construct local solutions.  The simplest version is for
uniformly small local oscillations as in the above theorem.  This result can be extended to uniform manifolds: We only need uniform local
coordinate maps. The uniqueness argument is elementary but delicate.

The estimates mentioned at the end of Section \ref{section2} required Calderon-Zygmund type estimates and pointwise
bounds of the heat kernel. Both are available for boundary value
problems in a half space.

Now consider a parabolic equation in a bounded domain with smooth boundary.
Locally we can flatten the boundary. We take the half space problem as model,
and consider the bounded domain with smooth boundary as a uniform manifold.

This allows to deal with  Dirichlet boundary conditions and
conormal boundary conditions
\[  \sum_{i,j=1}^d \nu_i a^{ij}(t,x,u) \partial_j u = \sum_i g_i(t,x,u) \partial_i u
+ f(t,x,u) \]
where we assume (with  $\nu$ denoting the exterior normal vector)
\[ \sum_{i=1}^n g_i(t,x,u) \nu_i(x)  = 0, \]
which expresses that $(g_i)$ has values in the tangent space of the boundary.

\subsection{Systems}
The same arguments apply to systems of equations
\[ u^k_t - \sum_{i,j=1}^d \sum_l \partial_i a^{ij}_{kl}(t,x,u) \partial_j u^l = f^k (t,x,u,\nabla u) \]
as soon as the Calderon-Zygmund estimates and the Gaussian estimates
are available.

A sufficient condition is that
\[  \sum_{i,j,k,l}  a^{ij}_{kl} A^k_i A^l_j \ge \lambda^{-1}  |A|^2 \]
holds uniformly. This implies that
\begin{equation} \label{positivity}
 \int \sum_{i,j,k,l} \tilde a^{ij}_{kl} (x) \partial_i \phi^k \partial_j \phi^l dx \ge \lambda^{-1}  \Vert |\nabla \phi| \Vert_{L^2}^2
\end{equation}
where $\tilde a^{ij}_{kl}(x) = a^{ij}_{kl} (t,x, v(t,x)) $ and
$v(t,x)$ is the mean of the initial data on a ball of radius $\sqrt{t}$ and
$ 0  < t \le T$.

The positivity condition \eqref{positivity} implies rank 1 positivity,
\[ \sum_{i,j,k,l} a^{ij}_{kl} \eta^k \eta^l \xi_i \xi_j \ge \lambda^{-1} |\xi|^2 |\eta|^2. \]
 On the other hand, for uniformly continuous coefficients, rank-1 positivity implies
\begin{equation} \label{positivity2}
 \int \sum_{oi,j,k,l} \tilde a^{ij}_{kl} (x) \partial_i \phi^k \partial_j \phi^l dx
 \ge (2\lambda)^{-1} \Vert \nabla \phi \Vert_{L^2}^2 - C \Vert \phi
 \Vert_{L^2}^2,
\end{equation}
but for discontinuous coefficients no good algebraic characterization
of the coefficients satisfying \eqref{positivity2} 
seems to be available.

In any case \eqref{positivity} for $a^{ij}_{kl}(t,x) = a^{ij}_{kl}(t,x,w(x))$, 
where $w$ is the heat extension of the initial data, is enough  to
handle small $BMO$ perturbations of uniformly continuous initial data.
The situation is different for small $L^\infty$ perturbations. Here we
deal with small perturbations of uniformly continuous coefficients and
rank 1 positivity suffices.

There is an important and natural weaker notion of parabolicity.
Assume that the coefficients $a^{ij}_{kl}$ have values in a compact set $K$ of tensors.
We call the equation parabolic if in this compact set $K$, for all
$\xi\in\R^d\backslash \{0\}$ the matrix
\[ A_{kl}(\xi) = \sum_{i,j=1}^d a^{ij}_{kl} \xi_i \xi_j\]
has its spectrum in the open left complex half plane.

The situation of boundary value problems is considerably more complex.
Again positivity is sufficient, and for Dirichlet boundary conditions
this is again equivalent to rank one positivity. In general no good
characterization of positivity is known, but there are many important sufficient
conditions, see Simpson and Spector \cite{MR713118}.

Again parabolicity in the sense of Solonnikov \cite{Sol65}  is sufficient for an
analogue of Theorem \ref{wellposed} for initial boundary value
problems for systems.

\subsection{More derivatives}
 Consider the parabolic equation
\[ u_t -\sum_{i,j=1}^d  a^{ij}(t,x,u,\nabla u ) \partial^2_{ij} u = f(t,x,u,\nabla u ) \]
in $\mathbb{R}^d$ and let $T > 0$ be given. We assume boundedness with a parameter
$\varepsilon$,
\[  \Vert a^{ij} \Vert_{sup} \le \lambda , \]
\[ \Vert f \Vert_{sup} \le \varepsilon T^{-1/2} \]
parabolicity,
\[  \sum_{i,j=1}^d a^{ij}(t,x,u,p)  \xi_i\xi_j \ge |\xi|^2 / \lambda \]
and Lipschitz continuity,
\[  |a^{ij}(t,x,u,p) - a^{ij}(t,x,v,q)| \le L (  |p-q|+ T^{-\frac12}  |u-v|) \]
\[  |f(t,x,u,p) - f(t,x,v,q)| \le \varepsilon  (T^{-\frac12}|p-q|+ T^{-1}  |u-v|). \]
We assume again  locally small oscillation:
There exists $\delta$ depending only on  $\lambda$ with
\[  \sup_{|x-y|\le \sqrt{T}, 0\le t,s \le T} |a^{ij}(t,x,u,p) -a^{ij}(s,y,v,q)| \le \delta \]
and regularity with $k \ge 1$,
\begin{equation}
 T^{|\alpha|/2+ l/2} |\partial_x^\alpha \partial_u^l \partial_p^\beta  a^{ij} | \le c
\end{equation}
and, with $ k\ge 1$,
\begin{equation}
 T^{\frac12+|\alpha|/2+ l/2} |\partial_x^\alpha \partial_u^l \partial_p^\beta  f  | \le c
\end{equation}
for $|\alpha| \le k$.

\begin{theorem}
\label{wellposed2}  There exist $ \delta >0$ such that for all $L >0$ there is
  $\varepsilon_0>0$ so that, if  $T>0$, 
\[
|\nabla u_0(x)-\nabla u_0(y)| \le \varepsilon \le \varepsilon_0  \text{ for  } |x-y| \le \sqrt{T}
\]
and if the assumptions above are satisfied
then there is a unique continuous solution $u$ up to time $T$ which satisfies
\[  t^{-\frac12} |(t^{\frac12} \partial_x)^\alpha u(t,x) | \le c_{\alpha}  \varepsilon \]
for $1\le |\alpha| \le 1+k$. The solution is analytic with respect
to $x$ if $a^{ij}$  and $f$ are analytic.
If $a^{ij}$ and $f$ are analytic with respect to all variables then there exist $c$ and $R$ so that
\[   t^{-\frac12} |(t^{\frac12}  \partial_x)^\alpha ((t\partial_t)^j u(t,x) | \le c (|\alpha|+j)! R^{-(|\alpha|+j)}  \varepsilon \]
for $|\alpha|+ j \ge 1$,
where $c$ and $R$ are independent of $x,t,j$ and $\alpha$.
\end{theorem}

The mean curvature flow in arbitrary codimension provides an example of 
this structure. Here bounded first derivatives seem to be appropriate if one
wants to deal with the flow for graphs, see e.g. \cite{MR2012810}, \cite{MR2128604}. Note that Sverak \cite{zbMATH00036262} has constructed
Lipschitz continuous singular solutions to the stationary problem. Either
these solutions indicate that solutions to the parabolic equation become nonunique, or
that the smallness condition for the initial data is needed for solutions
the function space $X$.

Similarly we deal with the fully nonlinear equation
\[ u_t - F(t,x,u,\nabla u , \nabla^2 u ) = 0 \]
with initial data in $C^{1,1}$. We assume Lipschitz continuity
\[ |F(t,x,u,p,A)-F(t,x,u,p,B)| \le \lambda |A-B|  \]
and ellipticity
\[  F(t,x,u,p,A+B) - F(t,x,u,p,A) \ge \lambda^{-1} \lambda_{min}(B) \]
where $A$  is symmetric and  $B$ positiv definite with  $\lambda_{min}$ denoting the smallest eigenvalue. 
 The condition of locally small
oscillation takes the form:
\begin{equation}  
\begin{split} 
\sup\limits_{|s-t|\le T, |x-y|\le \sqrt{T}} |F(t,x,.,.,A+H) - F(t,x,.,.,A) 
\hspace{-4cm} & \\ & - (F(s,y,.,.,A+H)-F(s,y,.,.,A)) |  \le \delta |H| 
\end{split} 
\end{equation}
The Lipschitz condition involving a small parameter $\varepsilon$ is
\[ |F(t,x,u,p,A) - F(t,x,v,q,A)| \le \frac{\varepsilon}{T} |u-v| + \frac{\varepsilon}{\sqrt{T}} |p-q| \]
and
\[ |F(t,x,u,p,0) | \le \varepsilon.  \]
Let $k\ge 1$. The higher regularity condition is
\[ |(T^{\frac12} \partial_x)^\alpha (T \partial_u)^l (T \partial_p)^\beta (\partial_A)^\gamma F | \le L \]
for $l+ |\alpha|+|\beta|+ |\gamma| \le k$.

\begin{theorem}
\label{wellposed3}  There exists $\varepsilon_0>0$ so that, if
\[
|D^2 u_0(x)-D^2 u_0(y)| \le \varepsilon \le \varepsilon_0  \text{ for  } |x-y| \le \sqrt{T}
\]
then there is a unique continuous solution $u$ up to time $T$ which satisfies
\[  t^{-1}  |(t^{\frac12} \partial_x)^\alpha u(t,x) | \le c_{\alpha}  \varepsilon \]
for $2\le |\alpha| \le 2+k$. The solution is analytic with respect
to $x$ if $a^{ij}$  and $F$ are analytic.
If $a^{ij}$ and $F$ are analytic with respect to all variables then there exist $c$ and $R$ so that
\[  t^{-1} |(t^{\frac12}  \partial_x)^\alpha ((t\partial_t)^j u(t,x) | \le c (|\alpha|+j)! R^{-(|\alpha|+j)}  \varepsilon \]
where $c$ and $R$ are independent of $x,t,j$ and $\alpha$.
\end{theorem}

It is not clear whether the smallness condition is needed. Note
however that Nadirashvili and {Vl\u adu\c t} \cite{zbMATH05247547} and
{Nadirashvili}, {Tkachev} and {Vl\u{a}du\c{t}} \cite{zbMATH06094090}
  have constructed singular solutions in $C^{1,1}$. So again, either
  the parabolic flow is nonunique for large $C^{1,1}$ initial data, or
  the smallness assumption is needed.

\section{Applications}
\subsection{Navier-Stokes equations}
Consider the Navier-Stokes equations
\[
\begin{split}
u_t-\Delta u + u\nabla u + \nabla p = &  0 \\
 \nabla \cdot u = & 0.
\end{split}
\]
with divergence free initial data $u_0$.

Let $v$ be the caloric extension, i.e. the solution to the heat equation
of the initial data $u_0$. The Carleson measure characterization of the BMO norm is
\[ \Vert u_0 \Vert_{BMO} \sim \sup_{R,x} \left( R^{-d} \int_{0}^{R^2} \int_{B_R(x)} |\nabla v(t,y)  |^2 dy dt \right)^{\frac12} \]
which we use to define the local $BMO^{-1}$ norm by
\[ \Vert u_0 \Vert_{BMO^{-1}_T} \sim \sup_{R^2\le T ,x} \left( R^{-d} \int_{0}^{R^2} \int_{B_R(x)} |v(t,y)  |^2 dy dt \right)^{\frac12}. \]
 We also let
\[ \Vert u \Vert_{X_T}=  \Vert t^{\frac12}|u(t,x)|\Vert_{sup}  +\sup_{x,R\le \sqrt{T}} \left( R^{-d} \int_0^{R^2} \int_{B_R(x)}  |u|^2 dy dt \right)^{\frac12}.   \]

\begin{theorem} There exists $\varepsilon>0$ depending only on the space dimension $d$ so that given $u_0$ with $\Vert u_0 \Vert_{BMO^{-1}_T} < \varepsilon$ there exists and  a unique solution $u \in X_T$ up to time $T$ with
\[ \Vert u \Vert_{X_T} \le c \Vert u_0 \Vert_{BMO^{-1}}. \]
\end{theorem}

The solution is a classical solution for $T>0$. It assumes the initial data
in the weak sense. See Koch and Tataru \cite{MR1808843} for more details.

\subsection{Hamilton-Jacobi equations and harmonic map heat flow}
Consider
\[ u_t - \sum_{i,j=1}^d \partial_i a^{ij}(x,u) \partial_j u = \sum_{i,j=1}^d f^{ij}(u) \partial_i u \partial_j u  \]
on a bounded domain $\Omega$ with smooth boundary and homogeneous
Dirichlet initial data where the coefficients $a^{ij}$ are bounded,
uniformly elliptic, with uniformly bounded derivatives. Also $f$ is
supposed to be bounded with uniformly bounded derivatives. The
harmonic map heat flow is a particular example, for which the coefficients
$a^{ij}$ are independent of $u$. In this form the type of the
equations does not change when we change dependent and
independent variables.

\begin{theorem}
There exists $\varepsilon$ such that the following is true:
Let $\phi_0 \in C(\overline\Omega) $ satisfy $\phi_0=0$ at the boundary.
There exists $T>0$ such that whenever
\[ \Vert u_0 -\phi_0 \Vert_{BMO} \le \varepsilon \]
then there is a unique smooth solution up to time $T$.
\end{theorem}

Here we use the heat extension with Dirichlet boundary conditions to
define the space $BMO$. It is remarkable that the initial data is not
required to satisfy the boundary condition.

Let us consider an example on $B_1(0) \subset\R^2$: We want to solve the equation
\[ u_t - \Delta u = |\nabla u|^2  \]
with initial data
\[ u_0 (x)= \ln ( 1- \ln(|x|) )  \]
which is in $BMO$.

Our results yield the existence of a unique smooth solution which assumes the initial data in a weak sense. It is remarkable that the constant map $u(t,x)= u_0$ is also a weak solution.

The harmonic map heat flow on $(0,T) \times \R^d$ has been considered previously by
the authors \cite{MR2916362}, and with small BMO initial data by Wang \cite{MR2781584}. We extend these results to uniform manifolds.

\subsection{Ricci-DeTurck flow}
The Ricci flow
\begin{align}
\partial_t g=&-2\text{Ric}(g)\ \ \ \text{in}\ \ M^n\times (0,T)\ \ \ \text{and}\nonumber\\
g(0,\cdot)=&g_0,\label{Ricci}
\end{align}
is the most natural parabolic deformation of a metric on a Riemannian manifold. Due to the
invariance under coordinate changes it is not parabolic. DeTurck \cite{MR697987} introduced
a condition fixing the coordinates: He considered a Ricci flow coupled with 
the harmonic map heat flow with respect to a background metric. 
In local coordinates the Ricci-DeTurck flow  can be written as
\[
\begin{split}
(\partial_t -\nabla_a g^{ab} \nabla_b) g_{ij} =& -\nabla_a g^{ab}  \nabla_b g_{ij} - g^{kl}g_{ip}h^{pq}R_{jkql}(h)-g^{kl}g_{jp}h^{pq}R_{ikql}(h)\\
&\hspace{-3cm} +\frac{1}{2}g^{ab} g^{pq} \times \Big( \nabla_i g_{pa} \nabla_j g_{qb}+2\nabla_a g_{jp} \nabla_q g_{ib}-2\nabla_a g_{jp} \nabla_b g_{iq} \\
&\hspace{-3cm} \qquad -2\nabla_j g_{pa} \nabla_b g_{iq}-2\nabla_i g_{pa} \nabla_b g_{jq}\Big)
\end{split}
\]
where we use a fixed background metric $h$. This is a particular instance
of Theorem \ref{wellposed} when we require that the initial metric lies in a
compact convex set of positiv definite matrices. 

By Whitney's result \cite{MR1503315} we may approximate a uniform
$C^1$ Riemannian manifold by a uniform $C^k$ Riemannian
manifold. Altogether, using Theorem \ref{wellposed} we arrive at
\begin{theorem} Let $(M,g_0)$ be a uniform $C^1$ manifold with a uniformly continuous metric $g_0$. Choose an atlas which makes $M$ a  uniform $C^3$ manifold
with $h$  a $C^2$ background metric with uniformly bounded second derivatives.
Then there exist $\varepsilon>0$ (independent of $g_0$), $T>0$ and a 
continuous solution $g$ of the Ricci-DeTurck flow on $(0,T)\times M$ with $g(0,\cdot)=g_0$ and which satisfies
\[
t^{1/2} \| \nabla (g(t)-h)\|_{L^\infty} \le \varepsilon.
\]
Moreover the solution is unique among all other solutions satisfying
the same bound for the gradient.
\end{theorem}    
We note that there are several interesting existence results for the
Ricci flow under various curvature assumptions using more geometric
arguments by Cabezas-Rivas and Wilking \cite{arXiv11070606} and Simon
\cite{MR1957662}, \cite{MR2876261}.

Uniqueness results were previously obtained under some curvature
bounds by Chen and Zhu \cite{MR2260930}, Chen \cite{MR2520796} and
Kotschwar \cite{arxiv12063225}. 

\subsection{Asymptotics for fast diffusion}

Consider the fast diffusion equation
\[ u_t =\frac1m \Delta u^m \]
with $m <1$. Let
\begin{equation} \label{beta}  \beta =  (2-(1-m)d)^{-1} \end{equation} 
and
\[ u_B = (B+|x|^2)^{-\frac{1}{1-m}}. \]
Then
\[   u(t,x)  = t^{-\beta d} (B+\frac{|x|}{t^\beta}^2)^{-\frac{1}{1-m}} \]
is the Barenblatt solution.

Conformal coordinates lead to the equation
\[\begin{split}   v_t = &  \frac1m (B+|x|^2) \Delta v^m +\frac2{1-m} x\cdot
\nabla (v-2v^m) \\ & + \left(d +2\frac{B+|x|^2}{1-m}|x|^2 \right) (v-v^m)
\end{split}
\]
This equation is uniformly parabolic on the cigar manifold given by the Riemannian metric
\[ \delta_{ij} (B+|x|^2)^{-1} \]
 provided the relative size $v= u/u_B$ is bounded from below and
 above. It has been shown by Vazquez that under weak assumptions on the
 initial data $v\to 0$ uniformly in $x$ as $t \to \infty$. It is
 remarkable that the spectrum and the eigenfunctions of the linearization can
 be computed explicitly, see  Denzler and  McCann \cite{zbMATH02155995}.

Using the formulation on the manifold above but not the approach
discussed here, Denzler, McCann and the first author \cite{arxiv12046434} derived precise
information on the large time asymptotics from the information on the
linearized operator. Due to the fact that the cigar is noncompact there
are important issues about the continuous spectrum for which we refer the reader to \cite{arxiv12046434}.

\subsection{Perturbed traveling wave solutions to the porous medium equation}
The porous medium equation
\[ \rho_t = \Delta \rho^m \]
with $m>1$ is an idealized model for the propagation of
gas in a porous medium. It has special solutions: The Barenblatt solution
\[ \rho(t,x)  = t^{-\beta d} \Big(B-\frac{|x|^2}{t^{\beta}}\Big)_+^{\frac1{m-1}} \]
which has compact support in $x$ for fixed $t$. Here $\beta$ is defined by \eqref{beta}.

A second explicit solution is given by the traveling wave solution
\[ \rho(t,x)^{m-1} = ( t+x_n)_+. \]

The quantity
\[ v = \frac{m}{m-1} \rho^{m-1} \]
corresponds to the physical pressure. It satisfies formally
\[ v_t - (m-1) v \Delta v = |\nabla v|^2. \]

\begin{theorem}[Kienzler 2013]
 Suppose that the nonnegative function $\rho_0 : \R^d \to \R$
 satisfies
\[
 \left| \nabla \left( \frac{m}{m-1} \rho_0^{m-1}\right) - e_n \right|
 < \delta\]
on the set of positivity. Then the unique solution to the
 porous medium equation satisfies
 \[
 \left|   \nabla \left( \frac{m}{m-1} \rho^{m-1}\right) - e_n \right| < C\delta
\]
and
\[   t^{k+|\alpha| -1} \left| \partial_t^k \partial_x^{\alpha}  \rho^{m-1}\right|
\le c_{k+|\alpha|} \delta \]
where $\rho$ is positiv 
whenever $1\le |\alpha|\le 2$.
\end{theorem}

Existence and uniqueness of solutions to nonnegative initial data is
well understood with the final contribution of Dahlberg and Kenig. The
regularity of solutions is more difficult. There are local regular
solutions to regular initial data satisfying a suitable nondegeneracy
condition (see Daskalopoulos and Hamilton \cite{zbMATH01198486} and
Daskalopoulos, Hamilton and Lee \cite{zbMATH01820822}).

 The Aronson-Graveleau solutions \cite{zbMATH00165220} describe the
 selfsimilar filling of a hole by gas.  It is a consequence that at
 the time of the filling the pressure does not remain Lipschitz
 continuous.

 Describing the graph is equivalent to describing
the function. We describe the graph of $p$ as a graph of a function $v$ with
\[ x_n = p, \qquad y_n = w. \]
It is defined on the halfplane $x_n >0$. The traveling wave solution becomes
\[ y_n-t  \]
and
\[ v = w-(y_n-t) \]
satisfies   with
\[ \sigma = \frac{m-2}{m-1} >-1 \]
\[ \frac1{m-1} v_t - (x_n^{-\sigma } \sum_{j=1}^{d-1} \partial_j (x_n^{1+\sigma}  \partial_j v))
- x_n^{-\sigma} \partial_n \left( x_n^{1+\sigma } \frac{\partial_n v - \sum_{j=1}^{d-1} (\partial_j v)^2}{1+\partial_n v} \right) = 0
\]
in the upper half plane $x_n \ge 0 $. The result in  transformed coordinates reads as
\begin{theorem}[C. Kienzler]
There exists $\delta>0$ such that the following is true. Suppose that
\[ v_0 : H \to \R \]
satisfies
\[  |v_0(x)-v_0(y)| \le \varepsilon |x-y|. \]
 Then there is a unique solution which satisfies
\[    |t^{j+|\alpha|-1} \partial_t^j \partial_x^{\alpha}  v|\le  c \varepsilon  \]
whenever $1\le |\alpha|\le 2$.
\end{theorem}

For the proof we observe that the second order part of the operator 
\[    x_n^{-\sigma} \nabla ( x_n^{1+\sigma} \nabla u) \] 
is the second order part of the Laplace-Beltrami
operator on the upper half plane with the
Riemannian metric
\[   \langle u,v\rangle_{x} = x_n^{-1} u\cdot v.  \]
This is half way  between Euclidean space and the Poincar\'e half  plane.

On an abstract level the steps are the same as on $\mathbb{R}^d$.
\begin{enumerate}
\item The intrinsic geometry defines balls and space time cylinders. On $L^2( x_n^{\sigma})$ we obtain a self adjoint
semigroup.
\item Energy arguments give $L^2$ estimates with Gaussian weights, the Davies-Gaffney estimates for the analogue of the
heat equation
\[ (m-1) v_t - x_n \Delta v - (1+ \sigma) v_n = 0. \]
\item A local regularity gives pointwise bounds of derivatives for solutions to the homogeneous equation in cylinders.
\item Both together imply Gaussian estimates for the fundamental solution and its derivatives in the intrinsic geometry.
\item The Gaussian estimates and the energy estimates are good enough for the Calderon-Zygmund theory on spaces of homogeneous type.
\end{enumerate}

See \cite{kienzler} for a complete proof.

\subsection{Flat solutions to the  thin film equation}
Nonnegative solutions to the thin film equation
\[ h_t + \nabla (h \nabla \Delta h) =0 \]
supposedly describe the dynamics of thin films. While existence of
weak solutions is reasonably well understood there are only few
instances where uniqueness or higher regularity are known. This is
a question with relevance for modeling: The equation has solutions with
zero contact angle, and nonzero contact angle, and hence at least the
contact angle is needed for a complete description. Here we study
existence and uniqueness of a class of solutions with zero contact
angle.  The only previous uniqueness result with a moving contact line
is in this setting in one space dimension by Giacomelli, Kn\"upfer and
Otto \cite{zbMATH05344296}.

There is a trivial stationary solution
\[ h = ((x_n)_+)^2  \]
and we want to study solutions in a neighborhood of $h$.

\begin{theorem}[D. John]
Suppose that
\[ |\nabla \sqrt{h_0}-e_n | \le \delta. \]
Then there exists a unique solution $h$ which satisfies
\[ |\nabla \sqrt{h} -e_n | \le c \delta \]
and, for $1\le |\alpha|\le 2$
\[ t^{2k+|\alpha|-1} \left| \partial_t^k \partial_x^{\alpha} \sqrt{h} \right|
\le c (k,\alpha) \Vert \nabla \sqrt{h_0} -e_n \Vert_{sup}. \]
\end{theorem}

This formulation is slightly different from what is proven by D. John in his thesis \cite{john}, but his proof gives also the simpler statement above.
Again we transform the problem to a degenerate quasilinear problem on the upper half plane.

Let $\tilde h = h^{\frac12}$ and note that it solves the equation
\[\begin{split}  \partial_t \tilde h + & \tilde h^2 \Delta^2 \tilde h + 6h\nabla \tilde h \nabla \Delta \tilde h
 +
\tilde h (\Delta \tilde h)^2 + 2h |\Delta' \tilde h|^2 + 2 |\nabla \tilde h|^2\Delta \tilde h
\\  + &
4 \partial_i \tilde h \partial_j \tilde h \partial^2_{ij} \tilde h = 0.
\end{split}
\]

Letting
\[ w= y_n\qquad  x_n = \tilde h \]
we obtain with $u= w-x_n$
\[ u_t + L_0 u =  f_0[u] + x_n f_1[u]+ x_n^2 f_2[u] \]
where
\[L_0 = x_n^{-1} \Delta x_n^3 \Delta -4 \Delta_{R^{n-1}} .  \]

The abstract procedure is the same as for the porous medium equation, but filling in the details is demanding.

\end{document}